\documentclass{amsart}

\usepackage{fullpage, setspace}
\title{Equations resolving a conjecture \\ of Rado on partition regularity}
\date{\today}

\author{Boris Alexeev}
\address{Department of Mathematics\\
Princeton University\\
Fine Hall, Washington Road\\
Princeton, NJ 08544-1000}
\email{balexeev@math.princeton.edu}

\author{Jacob Tsimerman}
\email{jtsimerm@math.princeton.edu}

\theoremstyle{plain}  % default
\newtheorem {theorem}            {Theorem}

\theoremstyle{definition}
\newtheorem*{definition*}        {Definition}

\newtheorem*{conjecture*}        {Conjecture}
\theoremstyle{remark}
\newtheorem*{remark*}            {Remark}

\DeclareMathOperator{\ord}{ord}

\begin{document}
\begin{abstract}
  A linear equation $L$ is called \emph{$k$-regular} if every
  $k$-coloring of the positive integers contains a monochromatic
  solution to $L$.  Richard Rado conjectured that for every positive
  integer $k$, there exists a linear equation that is $(k-1)$-regular
  but not $k$-regular.  We prove this conjecture by showing that the
  equation $\sum_{i=1}^{k-1} \frac{2^i}{2^i-1} x_i = \left( -1 +
  \sum_{i=1}^{k-1} \frac{2^i}{2^i-1} \right ) x_0$ has this property.

  This conjecture is part of problem E14 in Richard~K. Guy's book
  ``Unsolved problems in number theory'', where it is attributed to
  Rado's 1933 thesis, ``Studien zur Kombinatorik''.
\end{abstract}
\maketitle

In 1916, Schur~\cite{Schur} proved that in any coloring of the
positive integers with finitely many colors, there is a monochromatic
solution to $x+y=z$.  In 1927, van der Waerden~\cite{vdW} proved his
celebrated theorem that every finite coloring of the positive integers
contains arbitrarily long monochromatic arithmetic progressions.  In
his famous 1933 thesis, Richard Rado~\cite{Rado} generalized these
results by classifying the systems of linear equations with
monochromatic solutions in every finite coloring.  His thesis also
contained conjectures regarding equations that \emph{do} have a finite
coloring with no monochromatic solutions.

\begin{definition*}
  A linear equation $L$ is \emph{$k$-regular} if every $k$-coloring of
  the positive integers contains a monochromatic solution to $L$.
\end{definition*}
\begin{remark*}
  Some authors require that the values of the variables be distinct in
  solutions to $L$.  We follow Rado and Guy in not including this
  extra condition.
\end{remark*}

\begin{conjecture*}[Rado, \cite{Rado} via \cite{Guy}\footnote{These
      authors could not verify whether the conjecture is present in Rado's
      thesis.}]
  For every positive integer $k$, there exists a linear equation that is
  $(k-1)$-regular but not $k$-regular.  In other words, $k$ is the least
  number of colors in a coloring of the positive integers with no
  monochromatic solution to $L$.
\end{conjecture*}

Fox and Radoi\v{c}i\'{c}~\cite{FR} conjectured that the family of
equations $M_k$ given by $\sum_{i=0}^{k-2} 2^i x_i = 2^{k-1} x_{k-1}$
has this property.  We prove Rado's conjecture by using the related
family of equations,
\[ \sum_{i=1}^{k-1}
  \frac{2^i}{2^i-1} x_i = \left( -1 + \sum_{i=1}^{k-1}
  \frac{2^i}{2^i-1} \right ) x_0,\]
which we denote by $L_k$.

\begin{theorem}
  The equation $L_k$ is $(k-1)$-regular but not $k$-regular.
\end{theorem}
\begin{remark*}
  This result and its proof carry over to colorings of the nonzero
  rationals.
\end{remark*}
\begin{proof}
  We first use the power of $2$-adic valuations to prove that there
  exists a $k$-coloring with no monochromatic solutions to $L_k$.  The
  idea of using valuations was proposed by Fox and Radoi\v{c}i\'{c} for
  the equation $M_k$; later, Alexeev, Fox, and Graham~\cite{AFG} proved
  that these colorings were actually minimal if $k \le 7$, but only
  conjectured the result in general.

  If $r$ is a nonzero rational number, let $\ord_2(r)$ denote the
  $2$-adic valuation of $r$, that is, the unique integer $m$ such that
  $r = 2^m \frac{a}{b}$ for odd integers $a$ and $b$; also, let
  $\ord_2(0)= \infty$.  Recall that $\ord_2$ satisfies the following two
  properties:
  \begin{enumerate}
  \item
    $\ord_2(xy) = \ord_2(x) + \ord_2(y)$,
  \item
    $\ord_2(x+y) \ge \min( \ord_2(x), \ord_2(y) )$, with equality if
    $\ord_2(x)$ does not equal $\ord_2(y)$.
  \end{enumerate}
  The latter is the (non-Archimedean) $2$-adic triangle inequality.

  Note that the $2$-adic valuation of the coefficient of $x_i$ in $L_k$
  is $i$.  This is immediate if $i>0$ while in the case of $i=0$, this
  follows from the $2$-adic triangle inequality because $-1$ has
  valuation $0$ while the rest of the terms in the summation have
  positive valuation.  We claim this implies that the $k$-coloring
  $C_k(r) = \ord_2(r) \bmod k$ has no monochromatic solutions to $L_k$.

  Assume to the contrary that $C_k$ admits a monochromatic solution to
  $L_k$.  Then the terms of $L_k$ would be a collection of numbers
  with distinct $2$-adic valuations that sum to zero, which is
  impossible by the $2$-adic triangle inequality.

  We now show there is no coloring to $L_k$ with fewer than $k$ colors.
  Indeed, in any coloring with no monochromatic solutions to $L_k$, the
  color of $x$ and $2^jx$ must be different if $0 < j < k$.  (In the
  referenced literature, the number $2^j$ is thus said to be a
  \emph{forbidden ratio}.)  If there were an $x$ so that $x$ and $2^j x$
  were the same color, then
  \[ x_i = \begin{cases}
    x & \text{if }i=j, \\
    2^j x & \text{if }i \ne j,
  \end{cases}
  \]
  would be a monochromatic solution to $L_k$; all but the
  $j^\text{th}$ terms cancel, leaving $\frac{2^j}{2^j-1} x = \left(-1
  + \frac{2^j}{2^j-1} \right) 2^j x$.  This implies that the $k$
  numbers $1$, $2$, $4$, $\dotsc$, $2^{k-1}$ are colored with distinct
  colors.
\end{proof}

\section*{Acknowledgments}

The authors wish to thank Owen Biesel for helpful comments on the
exposition.

\bibliographystyle{amsalpha}
%\bibliography{rado}

\providecommand{\bysame}{\leavevmode\hbox to3em{\hrulefill}\thinspace}
\providecommand{\MR}{\relax\ifhmode\unskip\space\fi MR }
% \MRhref is called by the amsart/book/proc definition of \MR.
\providecommand{\MRhref}[2]{%
  \href{http://www.ams.org/mathscinet-getitem?mr=#1}{#2}
}
\providecommand{\href}[2]{#2}

\end{document}